\documentclass{llncs}
\usepackage{llncsdoc}

\usepackage{endfloat}
\usepackage{graphicx}
\usepackage{url}
\usepackage{algorithm}
\usepackage{algorithmic}
\usepackage{epstopdf}
\usepackage[usenames]{color}
\usepackage{rotating,afterpage}
\usepackage{multicol} 
\usepackage{amssymb}
\usepackage{latexsym}
\usepackage{amsmath}
\def\argmax{\mathop{\rm arg\,max}}

\begin{document}

\title{Auto-validating von Neumann Rejection Sampling from Small Phylogenetic Tree Spaces}
\titlerunning{Rigorous Phylogenetic Sampling}
\author{Raazesh Sainudiin\inst{1}  \and Thomas York\inst{2}}
\institute{Department of Statistics, 1 South Parks Road, University of Oxford,\\ Oxford OX1 3TG, U.K.
\email{sainudii@stats.ox.ac.uk}
\and Department of Biological Statistics and Computational Biology, Cornell University, Ithaca, U.S.A. \email{tly2@cornell.edu}}

\maketitle

\begin{abstract}
In phylogenetic inference one is interested in obtaining samples from the posterior distribution over the tree space on the basis of some observed DNA sequence data.  The challenge is to obtain samples from this target distribution without any knowledge of the normalizing constant.  One of the simplest sampling methods is the rejection sampler due to von Neumann.  Here we introduce an auto-validating version of the rejection sampler, via interval analysis, to rigorously draw samples from posterior distributions, based on homologous primate mitochondrial DNA, over small phylogenetic tree spaces.
\end{abstract}

\section{INTRODUCTION}

Obtaining samples from a density $p(\theta) \triangleq p^*(\theta)/N_p$, 
where $\theta \in \mathbf{\Theta}$ and $\mathbf{\Theta}$ is a compact Euclidean subset, {\it i.e.,} $\mathbf{\Theta} \subset \mathbb{R}^n$, without any knowledge of the normalizing constant $N_p \triangleq \int_{\mathbf{\Theta}}{p^*(\theta)} \, d \theta$, is a basic problem in statistical inference.  The usual Monte Carlo methods via conventional floating-point arithmetic are typically non-rigorous.  We will concentrate on the rejection sampler due to von Neumann \cite{Neumann1963} and its rigorous extension for application in phylogenetics.  The standard approaches to sampling from the posterior over phylogenies rely on Markov chain Monte Carlo (MCMC) methods.  Despite their asymptotic validity, it is nontrivial to guarantee that an MCMC algorithm has converged to stationarity \cite{Hobert2001}, and thus MCMC convergence diagnostics on phylogenetic tree spaces are heuristic \cite{Mossel2005}.  Thus, until now, no rigorous methodology has existed for perfectly sampling from the posterior distribution over phylogenetic tree spaces, even for $3$ or $4$ taxa.  Here, we solve this rigorous posterior sampling problem over small phylogenetic tree spaces.

After a brief introduction to the rejection sampler (RS) in Sect.~\ref{S:RSIntro}, 
an interval version of this sampler is formalized in Sect.~\ref{S:MRS}.  This sampler is referred to as the Moore rejection sampler (MRS) in honor of Ramon E.~Moore who was one of the influential founders of interval analysis \cite{Moore1967}.  
In Sect.~\ref{S:JCTQ}, we rigorously draw samples from the posterior over small tree spaces.  We conclude in Sect.~\ref{S:C}.  Section \ref{S:AppA} summarizes our notation and gives a brief introduction to interval analysis, a prerequisite to understanding MRS.  In Sect.~\ref{S:AppB}, Lemma \ref{lemma1} shows that MRS produces independent samples from the desired target density and Lemma \ref{lemma2} describes the asymptotics of the acceptance probability for a refining family of MRSs.  Unlike many conventional samplers, each sample produced by MRS is equivalent to a computer-assisted proof that it is drawn from the desired target, up to the pseudo-randomness of the underlying, deterministic, pseudo-random number generator.  An open source {\tt C++} class library for MRS is publicly available from {\tt www.stats.ox.ac.uk/\url{~sainudii}/codes} .  

\section{Rejection Sampler (RS)}\label{S:RSIntro}
Rejection sampling  \cite{Neumann1963} is a Monte Carlo method to draw independent 
samples from a target probability distribution $p(\theta) \triangleq p^*(\theta)/N_p$, 
where $\theta \in \mathbf{\Theta} \subset \mathbb{R}^n$.  Typically the target $p$ is any density that is 
absolutely continuous with respect to the Lebesgue measure.  In most cases of interest we 
can compute the target shape $p^*(\theta)$ for any $\theta \in \mathbf{\Theta}$, 
but the normalizing constant $N_p$ is unknown.  
The von Neumann RS can produce samples from $p$ according to Algorithm \ref{A:RS} when provided with (i) a proposal density $q(\theta) = q^*(\theta)/N_q$ from which independent samples can be drawn, $N_q \triangleq \int_{\mathbf{\Theta}}{q^*(\theta)} \, d \theta$ is known, and $q^*(\theta)$ is computable for any $\theta \in \mathbf{\Theta}$ and (ii) a constant $c$ defining the envelope function $f_q(\theta) \triangleq c q^*(\theta)$, such that, 
\begin{equation}
f_q(\theta) \triangleq c q^*(\theta) \geq p^*(\theta), \forall \, \theta \in \mathbf{\Theta} \enspace . \label{E:EnvCond}
\end{equation}

\begin{algorithm}
\caption{von Neumann RS}
\label{A:RS}
\begin{algorithmic}
\STATE {
{\it input:} 
(1) a target shape $p^*$,
(2) a proposal density $q$,
(3) an envelope function $f_q$ and
(4) an integer $TRIALS_{MAX}$
}
\STATE {\it output:} a sample from $U$ distributed according to $p$

\STATE {\it initialize:} $TRIALS \Leftarrow 0$, $SUCCESS \Leftarrow false$
\REPEAT
\STATE DRAW $T \sim q$ \COMMENT{draw a sample from the random variable $T$ with distribution $q$}
\STATE DRAW $H$ $\sim$ $Uniform[0,f_q(T)]$, where $f_q(T) \geq p^*(T)$
\IF {$H \leq p^*(T)$}
\STATE $U \Leftarrow T$, $SUCCESS \Leftarrow true$
\ENDIF
\STATE $TRIALS \Leftarrow TRIALS+1$
\UNTIL{$TRIALS < TRIALS_{MAX}$ or $SUCCESS = true$}
\end{algorithmic}
\end{algorithm}

$U$ generated by the above algorithm is distributed 
according to $p$ \cite{Williams2001}.  
Observe that the probability $\mathbf{A}^p_{f_q}$ that a point proposed according to $q$ gets accepted 
as an independent sample from $p$ through the envelope function $f_q$ is the ratio of the integrals 
\[
\mathbf{A}^p_{f_q} = \frac{N_p}{N_{f_q}} \triangleq 
\frac{ \int_{\mathbf{\Theta}} {p^*(\theta) \, d \theta} }{\int_{\mathbf{\Theta}} {f_q(\theta) \, d \theta} } \enspace ,
\]
and the probability distribution over the number of samples from $q$ to obtain one sample from $p$ is 
geometrically distributed with mean $1/\mathbf{A}^p_{f_q}$ \cite{Williams2001}.  

\section{Moore Rejection Sampler (MRS)}\label{S:MRS}

Moore rejection sampler (MRS) is an auto-validating rejection sampler (RS).  It can 
produce independent samples from any target shape $p^*$ that has a well-defined natural 
interval extension $P^*$ (Definition \ref{D:NIE}) over a compact domain $\mathbf{\Theta}$.  
MRS is said to be auto-validating because it automatically obtains a proposal $q$ that is easy 
to simulate from, and an envelope $f_q$ that is guaranteed to satisfy the envelope 
condition \eqref{E:EnvCond}.  In summary, the defining characteristics and notations of MRS are:
\[
\begin{array}{lcl}
\text{Compact domain} & \qquad & \mathbf{\Theta} = [\underline{\theta}, \overline{\theta}]\\
\text{Target shape} & \qquad & p^*(\theta) :\mathbf{\Theta} \rightarrow \mathbb{R}\\
\text{Target integral} & \qquad & N_p \triangleq \int_{\mathbf{\Theta}} {p^*(\theta) \, d \theta} \\ 
\text{Target density} & \qquad & p(\theta) \triangleq \frac{p^*(\theta)}{N_p} :\mathbf{\Theta} \rightarrow \mathbb{R}\\
\text{Interval extension of } p^* & \qquad & P^*(\Theta) : \mathbb{I}\mathbf{\Theta} \rightarrow \mathbb{IR}\\
\text{Proposal shape} & \qquad & q^*(\theta) :\mathbf{\Theta} \rightarrow \mathbb{R}\\
\text{Proposal integral} & \qquad & N_q \triangleq \int_{\mathbf{\Theta}} {q^*(\theta) \, d \theta} \\ 
\text{Proposal density} & \qquad & q(\theta) \triangleq \frac{q^*(\theta)}{N_q} :\mathbf{\Theta} \rightarrow \mathbb{R}\\
\text{Envelope function} & \qquad & f_q(\theta) = c q^*(\theta)\\
\text{Envelope integral} & \qquad & N_{f_q} \triangleq \int_{\mathbf{\Theta}} {f_q(\theta) \, d \theta}  = c N_q \\
\text{Acceptance probability} & \qquad & \mathbf{A}^p_{f_q} = \frac{N_p}{N_{f_q}}\\
\text{Partition of } \mathbf{\Theta} & \qquad & \mathfrak{T} \triangleq \{ \, \Theta^{(1)}, \Theta^{(2)}, ..., \Theta^{(|\mathfrak{T}|)}  \, \}.
\end{array}
\] 

If $p^* \in \mathfrak{E}$, the class of elementary functions (Definition \ref{D:ElemFunc}), its natural interval extension $P^*$ is well-defined on $\mathbf{\Theta}$ and $\mathfrak{T} \triangleq \{ \, \Theta^{(1)}, \Theta^{(2)}, ..., \Theta^{(|\mathfrak{T}|)}  \, \}$ be a finite 
partition of $\mathbf{\Theta}$, then by Theorem \ref{3.1.11} we can enclose $p^*(\Theta^{(i)})$, 
i.e., the range of $p^*$ over the $i$-th element of $\mathfrak{T}$, with the interval 
extension $P^*$ of $p^*$.
\begin{equation}\label{E:PartCont}
p^*(\Theta^{(i)}) \subseteq P^*(\Theta^{(i)}) \triangleq [\underline{P}^*(\Theta^{(i)}),\overline{P}^*(\Theta^{(i)})], \, \forall \, 
i \in \{1,2, ..., |\mathfrak{T}| \} \enspace . 
\end{equation}
For the given partition $\mathfrak{T}$ we can construct a partition-specific proposal $q^{\mathfrak{T}}(\theta)$ 
as a normalized simple function over $\mathbf{\Theta}$,
\begin{equation}\label{E:qRS}
q^{\mathfrak{T}}(\theta) = \left( N_{q^{\mathfrak{T}}} \right)^{-1} \, 
\sum_{i=1}^{|\mathfrak{T}|} { \overline{P}^*(\Theta^{(i)}) \, \mathbf{1}_{ \{ \theta \ \in \ \Theta^{(i)} \} } } \enspace ,
\end{equation}
with the normalizing constant 
$N_{q^{\mathfrak{T}}} \triangleq \sum_{i=1}^{|\mathfrak{T}|} \left( d(\Theta^{(i)}) \cdot \overline{P}^*(\Theta^{(i)}) \right)$, where, $d({\Theta}) = d([\underline{\theta}, \overline{\theta}]) =  
\overline{\theta} - \underline{\theta}$ is the {\em diameter} of ${\Theta}$.
The next ingredient $f_{q^{\mathfrak{T}}}(\theta)$ for our rejection sampler can simply be
\begin{equation}\label{E:fqRS}
f_{q^{\mathfrak{T}}}(\theta) = 
\sum_{i=1}^{|\mathfrak{T}|} { \overline{P}^*(\Theta^{(i)}) \, \mathbf{1}_{ \{ \theta \ \in \ \Theta^{(i)} \} } } \enspace .
\end{equation}
The necessary envelope condition \eqref{E:EnvCond} is satisfied by $f_{q^{\mathfrak{T}}}(\theta)$
because of  \eqref{E:PartCont}.  Now, we have all the ingredients to perform a more efficient 
partition-specific Moore rejection sampling.  Lemma \ref{lemma1} shows that if the target 
shape $p^*$ has a well-defined natural interval extension $P^*$,
and if $U$ is generated according to Algorithm \ref{A:RS}, 
and if the proposal density $q^{\mathfrak{T}}(\theta)$  and the envelope function 
$f_{q^{\mathfrak{T}}}(\theta)$ are given by  \eqref{E:qRS} and \eqref{E:fqRS}, respectively, 
then $U$ is distributed according to the target $p$.
Note that the above arguments as well as those in the proof of Lemma \ref{lemma1} naturally extend 
when $\mathbf{\Theta} \subset \mathbb{R}^n$ for $n>1$.  In the multivariate case, $\Theta^{(i)} \in \mathbb{IR}^n$ 
(Definition \ref{D:Boxes}) is a box.  Thus, we naturally replace the diameter of an interval by the {\em volume} of a box 
$v(\Theta^{(i)}) \triangleq  {\prod_{k=1}^{n}{d(\Theta_k^{i})} }$.  The envelopes and proposals are now simple
functions over a partition of the domain into boxes.  Analogous to the univariate case, the accepted samples 
are uniformly distributed in the region $S \subset \mathbb{R}^{n+1}$ `under' $p^*$ and `over' $\mathbf{\Theta}$.  
Hence their density is $p$ \cite{Williams2001}.

Next we bound the acceptance probability $\mathbf{A}^p_{f_{q^{\mathfrak{T}}}} \triangleq \mathbf{A}^p_{\mathfrak{T}}$ 
for this sampler.  Due to the linearity of the integral operator and \eqref{E:PartCont},
\[
\begin{array}{lcl}
N_p & \triangleq & \int_{\mathbf{\Theta}} {p^*(\theta) \, d \theta} \\
    &  = & \sum_{i=1}^{|\mathfrak{T}|} \int_{\Theta^{(i)}} {p^*(\theta) \, d \theta} \\
    &\in & \sum_{i=1}^{|\mathfrak{T}|} \left( d(\Theta^{(i)}) \cdot P^*(\Theta^{(i)}) \right) \\
    & =  & [ \ \, \sum_{i=1}^{|\mathfrak{T}|} \left( d(\Theta^{(i)}) \cdot \underline{P}^*(\Theta^{(i)}) \right) , \
    \sum_{i=1}^{|\mathfrak{T}|} \left( d(\Theta^{(i)}) \cdot \overline{P}^*(\Theta^{(i)}) \right) \, \ ] \enspace .
\end{array}
\]
Therefore, 
\[
\mathbf{A}^p_{\mathfrak{T}} = \frac{N_p}{N_{f_{q^{\mathfrak{T}}}}} = 
\frac{N_p}{\sum_{i=1}^{|\mathfrak{T}|} \left( d(\Theta^{(i)}) \cdot \overline{P}^*(\Theta^{(i)}) \right)} \geq 
\frac{\sum_{i=1}^{|\mathfrak{T}|} \left( d(\Theta^{(i)}) \cdot \underline{P}^*(\Theta^{(i)}) \right)}
{\sum_{i=1}^{|\mathfrak{T}|} \left( d(\Theta^{(i)}) \cdot \overline{P}^*(\Theta^{(i)}) \right)} \enspace .
\]

 If $p^* \in \mathfrak{E_L}$, the Lipschitz class of elementary functions (Definition \ref{D:LipElemFunc}), then we might expect the enclosure of $N_p$ to be proportional to the mesh $w \triangleq \max_{i \in \{1,\dots,\mathfrak{T}\}}{d(\Theta^{(i)})}$ of the partition $\mathfrak{T}$.  Lemma \ref{lemma2} shows that if $p^* \in \mathfrak{E_L}$ and $\mathfrak{U}_W$ is a uniform partition of $\mathbf{\Theta}$ into $W$ intervals, then the acceptance probability $\mathbf{A}^p_{\mathfrak{U}_W} = 1 - \mathcal{O} (1/W)$.
Thus, the acceptance probability approaches $1$ at a rate that is no slower than linearly with the mesh.  
We can gain geometric insight into the sampler from an example.  The dashed lines of a given shade, 
depicting a simple function in Fig.~\ref{Fi:refine}, is a partition-specific envelope function \eqref{E:fqRS} 
for the target shape $s^{*}(x) = -\sum_{k=1}^5{ k \, x \, \sin{(\frac{k(x-3)}{3})}}$ over the domain $\mathbf{\Theta} = [-10,6]$ and its normalization gives the corresponding proposal function \eqref{E:qRS}.  As the refinement of $\mathbf{\Theta}$ 
proceeds through uniform bisections, the partition size increases as $2^i$, $i=1,2,3,4$.  Each of 
the corresponding envelope functions in increasing shades of gray can be used to draw auto-validated samples 
from the target $s(x)$ over $\mathbf{\Theta}$.  Note how the acceptance probability (ratio of the area below the target shape to that below the envelope) increases with refinement.

We studied the efficiency of uniform partitions for their mathematical tractability.  In practice, we may further 
increase the acceptance probability for a given partition size by adaptively 
partitioning $\mathbf{\Theta}$.  In our context, adaptive means the possible exploitation of any current 
information about the target.  We can refine the current partition $\mathfrak{T}_{\alpha}$ and 
obtain a finer partition $\mathfrak{T}_{\alpha^{\prime}}$ with an additional box
by bisecting a box $\Theta^{(*)} \in \mathfrak{T}_{\alpha}$ along the side with the maximal diameter.  
There are several ways to choose a $\Theta^{(*)} \in \mathfrak{T}_{\alpha}$ for bisection.  When 
$\Theta^{(i)} \in \mathbb{IR}^n$ has volume 
$v(\Theta^{(i)})$, an optimal choice for $\Theta^{(*)} =  \argmax_{\Theta^{(i)} \in \mathfrak{T}_{\alpha}} 
{  \left( v(\Theta^{(i)}) \cdot d(P^*(\Theta^{(i)}) \right) }$.  Under this partitioning scheme, we employ a priority queue to conduct sequential refinements 
of $\mathbf{\Theta}$.  This approach avoids the exhaustive $\argmax$ computations to obtain
the $\Theta^{(*)}$ for bisection at each refinement step.  Once we have any partition $\mathfrak{T}$ of $\mathbf{\Theta}$, we can efficiently sample 
$\theta \sim q^{\mathfrak{T}}$ given by \eqref{E:qRS} in two steps.  First we sample a box 
$\Theta^{(i)} \in \mathfrak{T}$ according to the discrete distribution $t(\Theta^{(i)})$,
\begin{equation}
t(\Theta^{(i)}) =
 \frac{v(\Theta^{(i)}) \cdot \overline{P}^*(\Theta^{(i)})}
 {\sum_{i=1}^{|\mathfrak{T}|}{v(\Theta^{(i)}) \cdot \overline{P}^*(\Theta^{(i)})}},
  \   \Theta^{(i)} \in \mathfrak{T}, 
\end{equation}
and then we choose a $\theta \in \Theta^{(i)}$ uniformly at random.  Sampling from large discrete 
distributions (with million states or more) can be made faster by preprocessing the probabilities 
and saving the result in some convenient lookup table.  This basic idea \cite{Marsaglia1963} 
allows samples to be drawn rapidly.  We employ a more efficient preprocessing strategy \cite{Walker1977} 
that allows samples to be drawn in constant time even for very large discrete distributions as 
implemented in the GNU Scientific Library \cite{GSL}.  Thus, by means of priority queues and lookup
tables we can efficiently manage our adaptive partitioning of the domain for envelope construction, 
and rapidly draw samples from the proposal distribution.  We used the Mersenne Twister 
random number generator \cite{Matsumoto1998} in this paper.  Our sampler class builds 
on {\tt C-XSC 2.0}, a {\tt C++} class library for extended scientific computing using interval 
methods \cite{Hofschuster2004}.  All computations were done on a 2.8 GHz Pentium IV machine with
1GB RAM.  Having given theoretical and practical considerations to our 
Moore rejection sampler, we are ready to draw samples from various targets. 

\section{Auto-validating Independent Posterior Samples from Triplets and Quartets}\label{S:JCTQ}
  
Inferring the ancestral relationship among a set of species based on their DNA sequences is a basic problem in phylogenetics \cite{Semple2003}.  One can obtain the likelihood of a particular phylogenetic tree that relates the species of interest by superimposing a simple Markov model of DNA substitution due to Jukes and Cantor \cite{Jukes1969} on that tree.  The length of an edge (branch length) connecting two nodes (species) in the tree represents the amount of evolutionary time (divergence) between the two species.  The likelihood function over trees obtained through a post-order traversal (e.g.~\cite{Felsenstein1981}) has a natural interval extension over boxes of trees \cite{SainudiinPhD2005}.  This allows us to draw samples from the posterior distribution over
a compact box specified by our prior distribution on the tree space using our MRS.  We assume a uniform prior over the possible unrooted topologies and a uniform product prior over all branch lengths in the range $[10^{-10}, 10]$.  We consider two mitochondrial DNA data sets.
\subsection{Chimpanzee, Gorilla, Orangutan and Gibbon}
Our posterior distribution is based on the data from an $895$ bp long homologous segment of mitochondrial DNA from chimpanzee, gorilla, orangutan, and gibbon, containing the genes for three transfer RNAs and parts of two proteins \cite{Brown1982}.  Under the assumption of independence across sites, the sufficient statistics only comprise of the distinct site patterns and their counts.  The data for chimpanzee, gorilla and orangutan can be summarized by the following 29 distinct site patterns and counts: 
{\small
\begin{center}
\begin{verbatim}
site       :                   1 1 1 1 1 1 1 1 1 1 2 2 2 2 2 2 2 2 2 2
pattern    : 1 2 3 4 5 6 7 8 9 0 1 2 3 4 5 6 7 8 9 0 1 2 3 4 5 6 7 8 9
 . . . . . . . . . . . . . . . . . . . . . . . . . . . . . . . . . . .
chimpanzee : a g c t a t c a c c c a t c t g c c g t a c t a a g c g t
gorilla    : a g c t g t t a t c a a c a c g c a a a a t c c g g t a t
orangutan  : a g c t a c c g t t c c c a t a a t a a t a a a g c g c a
 . . . . . . . . . . . . . . . . . . . . . . . . . . . . . . . . . . .
site       : 2 7 2 1 1 3 1 1 9 2 1 8 2 3 1 8 7 1 9 2 4 2 1 2 1 1 2 1 3
pattern    : 3 1 2 6 3 1 6 8   0     2   0  
counts     : 2   9 8
\end{verbatim}
\end{center}}

\begin{figure}
\makebox{\centerline{\includegraphics{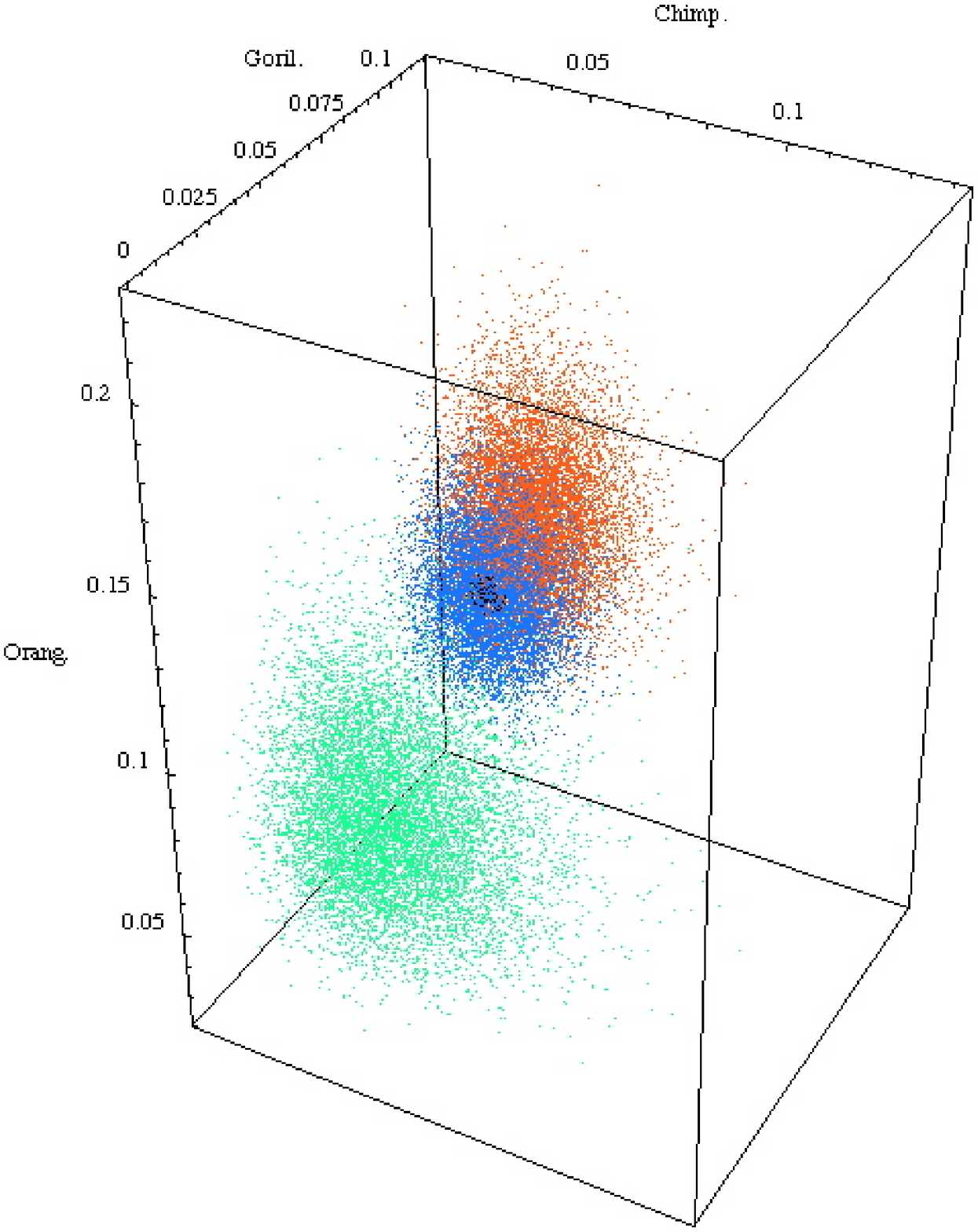}}} 
\caption[]{$10000$ Moore rejection samples from the posterior distribution over the three branch 
lengths of the unrooted phylogenetic tree space of Chimpanzee, Gorilla and Orangutan based on their homologous mitochondrial DNA sequence of length $895$ base pairs (blue dots), the tRNA-coding sequence with $198$ base pairs (green dots) and the protein-coding sequence with $697$ base pairs (red dots).  The verified maximum likelihood estimate is the large black dot within the blue dots.  \label{JC3MRS}}
\end{figure}

\begin{figure}
\makebox{\centerline{\includegraphics{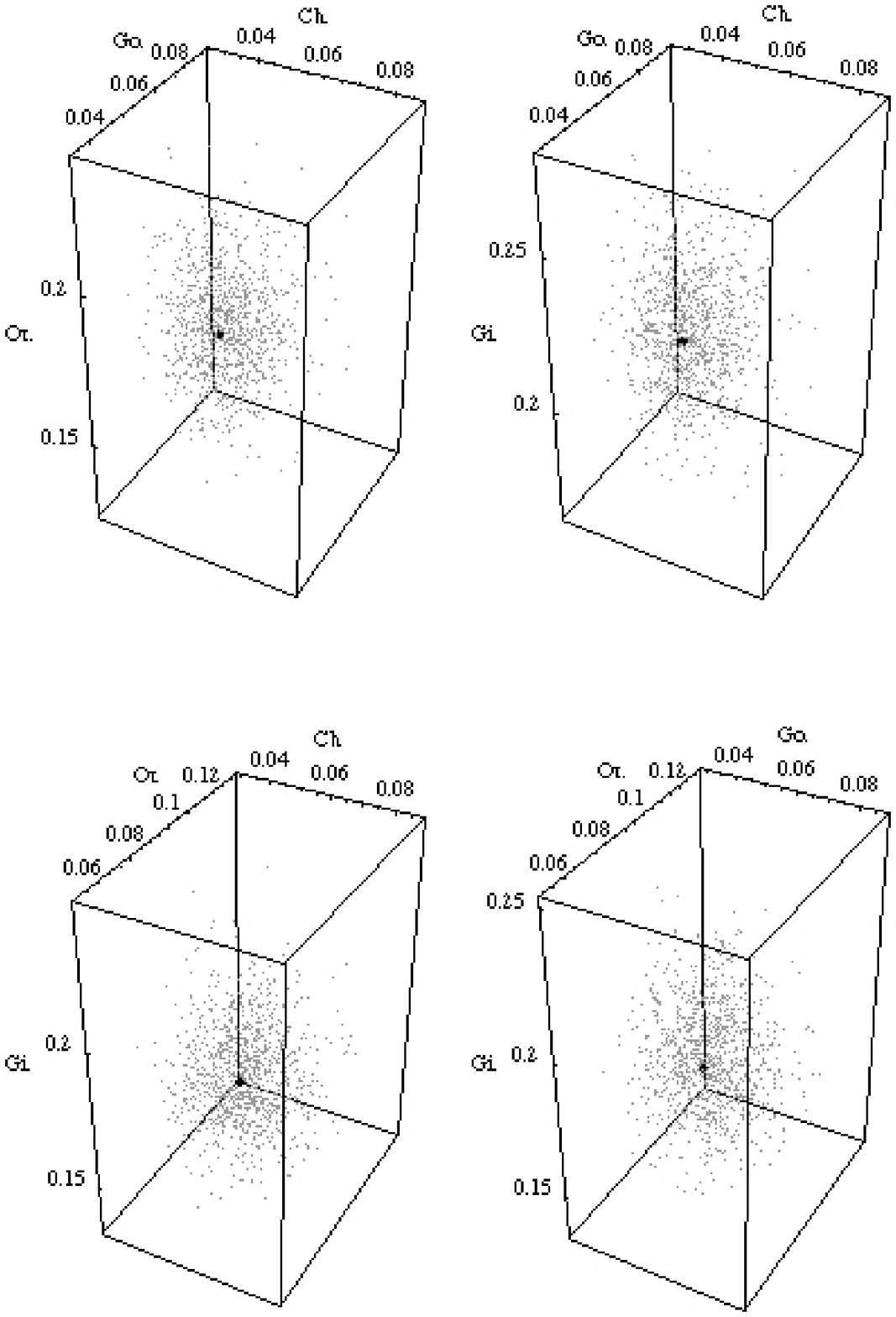}}} 
\caption[]{$1000$ Moore rejection samples (gray dots) from the posterior distribution over the unrooted quartet tree space of Chimpanzee (Ch.), Gorilla (Go.), Orangutan (Or.) and Gibbon (Gi.) depicted over the four marginal triplet branch lengths.  The verified marginal maximum likelihood estimates are the large black dots.  \label{JC4MRS}}
\end{figure}

In the above data set, the first column (1.aaa.232) expresses that there are $232$ site patterns with nucleotide `a' in all three species, ..., and the last column (29.tta.3) expresses that there are 3 site patterns with nucleotide `t' in chimpanzee and gorilla, and nucleotide `a' in orangutan.  The data for all four primates can be summarized by $61$ distinct site patterns as parsed in \cite{SainudiinPhD2005}.  $10000$ independent samples were drawn in $942$ CPU seconds from the posterior distribution over Jukes-Cantor triplets, i.e. unrooted trees with three edges corresponding to the three primates emanating from their common ancestor.  Figure \ref{JC3MRS} shows these samples (blue dots) scattered about the verified global maximum likelihood estimate (MLE) of the triplet obtained in \cite{SainudiinPhD2005} and subsequently confirmed algebraically in \cite{Hosten2005}.  We also drew $10000$ independent samples from the posterior based on the $198$ tRNA-coding DNA sites (green dots in Fig.~\ref{JC3MRS}) as well as from that based on the remaining $697$ protein-coding sites (red dots in Fig.~\ref{JC3MRS}).  The former posterior samples, corresponding to the tRNA-coding sites, are more dispersed than the posterior samples based on the entire sequence.  This is due to the smaller number of tRNA-coding sites making the posterior less concentrated.  We were able to reject the null hypothesis of homogeneity between the posterior samples based on the tRNA-coding sites and that based on the protein-coding sites at the $10\%$ significance level (P-value $= 0.06$ from a non-parametric bootstrap of Hotelling's trace statistic based on $100$ random permutations of the sites).  Any biological interpretation of this test must be done cautiously since the Jukes and Cantor model employed here forbids any transition$:$transversion bias that is reportedly relevant for this data \cite{Brown1982}.
We were able to draw samples from Jukes-Cantor quartets by adding the homologous sequence of the Gibbon.  Now, the problem is a more challenging because there are three distinct tree topologies in the unrooted, bifurcating, quartet tree space, and each of these topologies has five edges.  Thus, the domain of quartets is a piecewise Euclidean space that arises from a fusion of $3$ distinct five dimensional orthants.  Since the post-order
traversals specifying the likelihood function are topology-specific, we extended the likelihood over a compact box of quartets in a topology-specific manner.  The computational time was about a day and a half to draw $10000$ samples from the quartet target due to low acceptance probability of the naive likelihood function based on the $61$ distinct site patterns.  All the samples had the topology which grouped Chimp and Gorilla together, i.e. ((Chimp, Gorilla), (Orangutan, Gibbon)).  The samples were again scattered about the verified global MLE of the quartet \cite{SainudiinPhD2005}.  The marginal triplet trees (gray dots) within the sampled quartets are depicted in Fig.~\ref{JC4MRS}.  This quartet likelihood function has an elaborate DAG (Definition \ref{D:DAG}) with numerous operations.  When the data got compressed into sufficient statistics through algebraic statistical methods \cite{Casanellas2005}, the efficiency increased tremendously (e.g.~ for triplets the efficiency increases by a factor of $3.7$).  This is due to the number of leaf nodes in the target DAG, which encode the distinct site patterns of the observed data into the likelihood function, getting reduced from $29$ to $5$ for the triplet target and from $61$ to $15$ for the quartet target \cite{Casanellas2005}.  Poor sampler efficiency makes it currently impractical to sample from trees with five leaves and $15$ topologies (see Sect.~\ref{S:C} for a discussion on improvements).  However, one could use such triplets and quartets drawn from the posterior distribution to stochastically amalgamate and produce estimates of larger trees via fast amalgamating algorithms \cite{Strimmer1996,Levy2005}, which may then be used to combat the slow mixing in MCMC methods \cite{Mossel2005} by providing a good set of initial trees.

\subsection{Neandertal, Human and Chimpanzee}
We used the whole mitochondrial genome shotgun sequence (gi$| 115069275$) of a Neandertal fossil Vi-80, from Vindija cave, Croatia  \cite{Green2006}, and its homologous sequence in a human (gi$| 13273200$) and a chimpanzee (gi$| 1262390$), as summarized by the $15$ sufficient site patterns and their counts below, to conduct statistical inference about the human-neandertal divergence time.  

{\small
\begin{center}
\begin{verbatim}
site       :                   1 1 1 1 1 1 
pattern    : 1 2 3 4 5 6 7 8 9 0 1 2 3 4 5
 . . . . . . . . . . . . . . . . . . . . . 
neandertal : t t c a g g t g t c a a c a a
human      : t t c a g g t a c c a g t a g
chimpanzee : t c c a g a a a t t g a c t g
 . . . . . . . . . . . . . . . . . . . . . 
site       : 6 1 6 6 4 1 2 1 2 1 1 1 1 1 1
pattern    : 0 4 0 8 5 0       4 5
counts     : 5   3 5 0
\end{verbatim}
\end{center}}

\begin{figure}
\makebox{\centerline{\includegraphics[height=8.0cm]{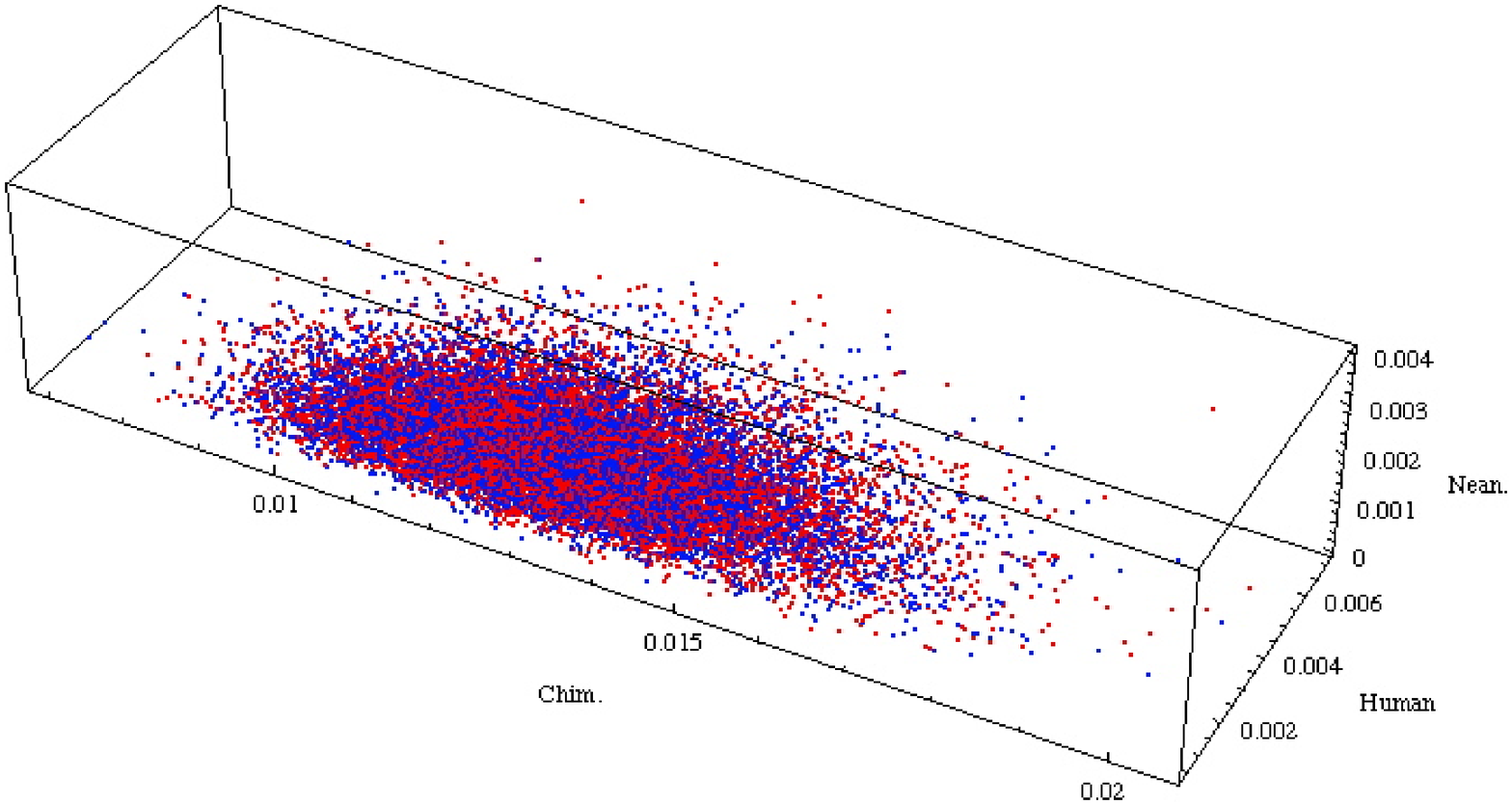}}} 
\caption[]{$10000$ Moore rejection samples each from the posterior distribution over the three branch lengths of the rooted phylogenetic tree space of Chimpanzee, Human and Neandertal under the Jukes-Cantor model (blue dots) and the HKY model (red dots) \label{JC3HKY3NHCRoot}}
\end{figure}

We drew $10000$ auto-validating independent samples from each of three posterior distributions; (1) over the space of unrooted triplets under the Jukes-Cantor model in $312$ CPU seconds, (2) over the clocked and rooted triplets under a Jukes-Cantor model in $375$ CPU seconds and (3) over the clocked and rooted triplets under a more general mutational model due to Hasegawa, Kishino and Yano (HKY) \cite{HKY85} in $1.2$ CPU hours.  In the HKY model we used the empirical nucleotide frequencies from the data ($\pi_{T} = 0.2588$,  $\pi_{C} = 0.2571$, $\pi_{A} = 0.2916$, $\pi_{G} =  0.1925$) and a hominid-specific transition/transversion rate of $2.0$.  Unlike the Jukes-Cantor model, all $15$ distinct site patterns are minimally sufficient under the HKY model and this is reflected in its longer CPU time.  Both models gave similar posterior samples over rooted triplets, as shown in Fig.~\ref{JC3HKY3NHCRoot}.

\begin{figure}
\makebox{\centerline{\includegraphics[height=8.0cm]{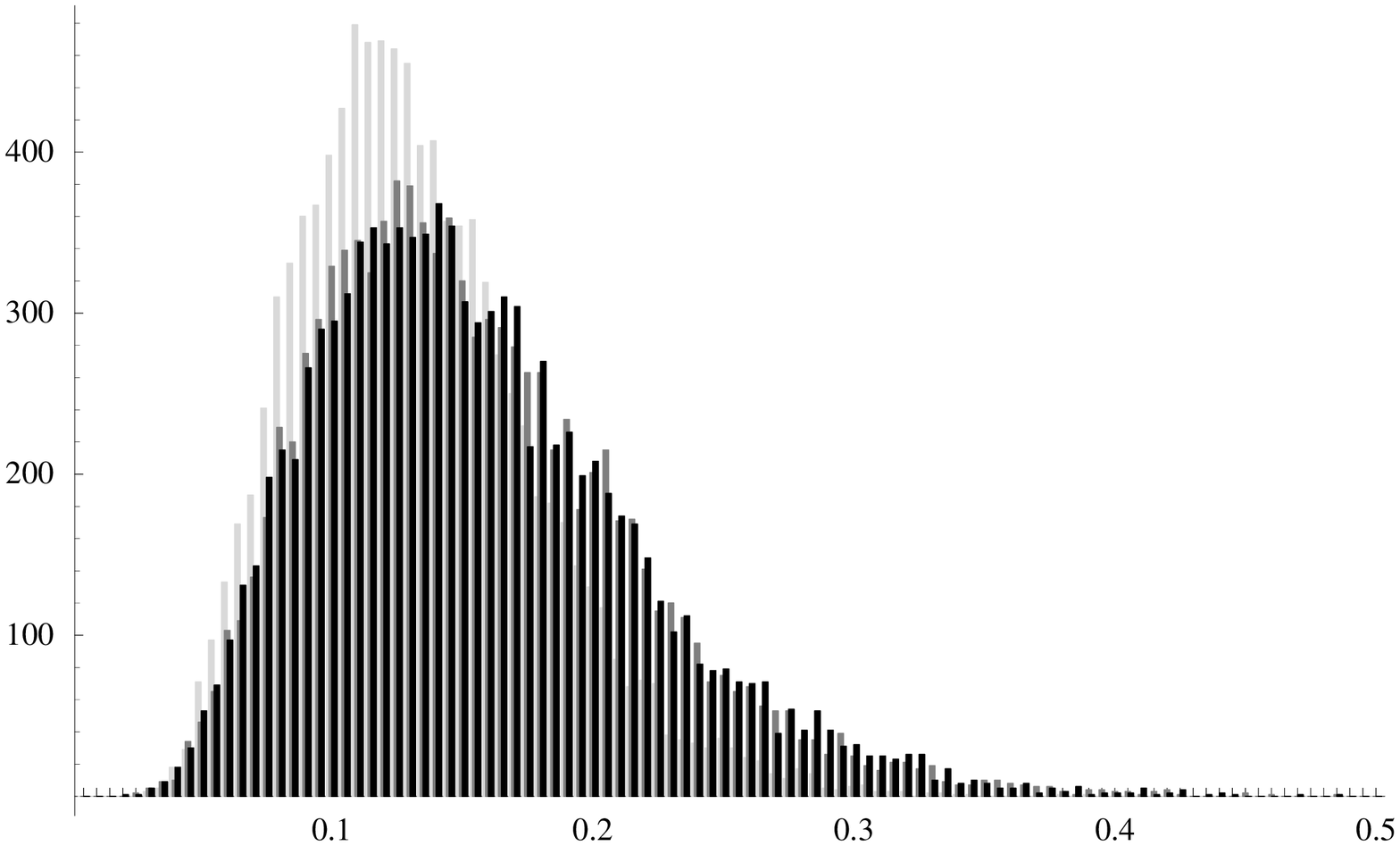}}} 
\caption[]{Posterior distribution over the human-neandertal divergence time relative to the human-chimp divergence time based on $10000$ independent samples from the (1) Midpoint-rooted tree estimates of the unrooted triplets under the Jukes-Cantor model (light gray), (2) rooted triplets under the Jukes-Cantor model (dark gray), and (3) rooted triplets under HKY model (black)  \label{HNDiv}}
\end{figure}

We transformed the three posterior distributions over the triplet spaces; (1) unrooted Jukes-Cantor triplets that were rooted using the mid-point rooting method, (2) rooted Jukes-Cantor triplets and (3) rooted HKY triplets, respectively, into three posterior distributions over the human-neandertal divergence time relative to the human-chimp divergence time (Fig.~\ref{HNDiv}).  The corresponding posterior quantiles ($\{ 5\%$, $50\%$ $, 95\% \}$) for the human-neandertal divergence times are $\{ 0.0643 $ $, 0.125 $ $, 0.214 \}$, $\{ 0.0694 $ $, 0.142 $ $, 0.263 \}$ and $\{ 0.0682 $ $, 0.143$ $, 0.268 \}$, respectively.  We constrained the neandertal lineage to be a fraction of the human lineage in branch length in order to estimate the age of the neandertal fossil from the rooted HKY triplets.   The posterior quantiles of the fossil date in units of human-chimp divergence is 
$\{ 0.00685$ $, 0.0666$ $, 0.195 \}$.  The estimate of $38,310$ years based on carbon-14 accelerator mass spectrometry \cite{Green2006} is within our $[5\%,95\%]$ posterior quantile interval for the fossil date, provided the human-chimp divergence estimates ranges in $[196103, 5.6 \times 10^6]$.  Thus, reasonable bounds for the human-chimp divergence are $4 \times 10^6$ and $5.6 \times 10^6$ years.  Based on these calendar year estimates, we transformed the posterior quantiles of the human-neandertal divergence times from the rooted HKY triplets into $\{ 272680$ $, 571124$ $, 1073375 \}$ and $\{ 381752$ $, 799574$ $, 1502724 \}$, respectively.  Our $[5\%,95\%]$ posterior intervals contain the interval estimate of $[461000,825000]$ years reported in \cite{Green2006}.  However, our confidence intervals are from perfectly independent samples from the posterior and account for the finite number of neandertal sites that were successfully sequenced, unlike those obtained on the basis of a bootstrap of site patterns \cite{Efron1996} or heuristic MCMC \cite{Hobert2001}.   Unfortunately, our human-neandertal divergence estimates are overestimates as they ignore the non-negligible time to coalescence of the human and neandertal homologs within the human-neandertal ancestral population.  Improvements to our estimates based on the other $310$ human and $4$ chimpanzee homologs reported in \cite{Green2006} may be possible with more sophisticated models of populations within a phylogeny and need further investigation.  

\section{Conclusion} \label{S:C}

Interval methods provide for a rigorous sampling from posterior target densities over small phylogenetic tree spaces.  When one substitutes conventional floating-point arithmetic for real arithmetic in a computer and uses discrete lattices to construct the envelope and/or proposal, it is generally not possible to guarantee the envelope property, and thereby ensure that samples are drawn from the desired target density, except in special cases \cite{GilksWild1992}.  Thus, the construction of the Moore rejection sampler through interval methods, that enclose the target shape over the entire real continuum in any box of the domain with machine-representable bounds, in a manner that rigorously accounts for all sources of numerical errors (see \cite{Kulisch2001} 
for a discussion on error control), naturally guarantees that the Moore rejection samples are independent draws from the desired target.  Moreover, the target is allowed to be multivariate and/or non-log-concave with possibly `pathological' behavior, as long as it has a well-defined interval extension.  

The efficiency of MRS is not immune to the curse of dimensionality and target DAG complexity.  
When the DAG for the likelihood gets large, its natural interval extension can have terrible over-enclosures of the true range, which in turn forces the adaptive refinement of the domain to be extremely fine for efficient envelope construction.  Thus, a naive application of interval methods to targets with large DAGs can be terribly inefficient.  In such cases, sampler efficiency rather than rigor is the issue.  Thus, one may fail to obtain samples in a reasonable time, rather than (as may happen with non-rigorous methods) produce samples from some unknown and undesired target. 
There are several ways in which efficiency can be improved for such 
cases.  First, the particular structure of the target DAG should be exploited to avoid any redundant computations.
For example, algebraic statistical methods can be used to find sufficient statistics to dissolve symmetries 
in the DAG as done in Sect.~\ref{S:JCTQ}.  Second, we can further improve efficiency by limiting ourselves to differentiable targets in $C^n$.  Tighter enclosures of the range $p^*(\Theta^{(i)})$ with $P^*(\Theta^{(i)})$ can come from the enclosures of Taylor expansions of $p^*$ around the midpoint $m(\Theta^{(i)})$ through interval-extended automatic differentiation (see \cite{Kulisch2001}) that can then yield tighter estimates of the integral enclosures \cite{Tucker2004}.  Third, we can employ pre-processing to improve efficiency.  For example, we can pre-enclose the range of a possibly rescaled $p^*$ over a partition of the domain and then obtain the enclosure of $P^*$ over some arbitrary $\Theta$ through a combination of hash access and hull operations on the pre-enclosures.  Such a pre-enclosing technique reduces not only the overestimation of target shapes with large DAGs but also the computational cost incurred while performing interval operations with processors that are optimized for floating-point arithmetic.  Fourth, efficiency at the possible cost of rigor can also be gained  (up to $30\%$ ) by foregoing directed rounding during envelope construction.  

\section{Acknowledgments}
This was supported by a joint NSF/NIGMS grant DMS-02-01037.  R.S.~is a Research Fellow of the Royal Commission for the Exhibition of 1851.  Many thanks to Rob Strawderman and Warwick Tucker for constructive comments on the sampler and Jo Felsenstein for clarifying the transition probabilities under the HKY model.

\section{Appendix A}\label{S:AppA}

\begin{definition}
Let $X \triangleq [\underline{x}, \overline{x}]$ be an interval in 
$\mathbb{IR} \triangleq \{ [ \underline{x}, \overline{x}] : \underline{x} \leq \overline{x}, \underline{x}, \overline{x} \in \mathbb{R} \}$
\end{definition}

\begin{definition}[Interval arithmetic]\label{Df:intarith}
If the binary operator $\star$ is one of the elementary arithmetic 
operations $\{ +,-,\cdot,/\}$, then we define an arithmetic on operands in 
$\mathbb{IR}$ by
\[
X \star Y \triangleq \{x \star y : x \in X , y \in Y\}
\]
with the exception that $X/Y$ is undefined if $0 \in Y$.
\end{definition}
\begin{theorem}\label{Th:IRarith}
Arithmetic on the pair $X,Y \in \mathbb{IR}$ is given by:
\[
\begin{array}{lcl}
X + Y  & = & [\underline{x} + \underline{y}, \overline{x} + \overline{y}] \\
X - Y & = & [\underline{x} - \overline{y}, \overline{x} - \underline{y}] \\
X \cdot Y & = & [\min \{\underline{x} \underline{y}, \underline{x} \overline{y},
\overline{x} \underline{y}, \overline{x} \overline{y}\},
\max \{\underline{x} \underline{y}, \underline{x} \overline{y},
\overline{x} \underline{y}, \overline{x} \overline{y}\}], \notag \\
X / Y & = & X \cdot [1/\overline{y}, 1/\underline{y}], \text{ provided, } 0 \notin Y. \notag
\end{array}
\]
\end{theorem}
{\bf Proof} (cf.~\cite{Tucker2004}):
Since any real arithmetic operation $x \star y$, where $\star \in \{ +,-,\cdot,/\}$ and
$x,y \in \mathbb{R}$, is a continuous function
$x \star y\triangleq\star(x,y): \mathbb{R} \times \mathbb{R} \rightarrow \mathbb{R}$,
except when $y=0$ under $/$ operation.  Since $X$ and $Y$ are simply connected compact 
intervals, so is their product $X \times Y$.  On such a domain $X \times Y$, the continuity 
of $\star(x,y)$ (except when $\star = /$ and $0 \in Y$) ensures the attainment of a minimum, 
a maximum and all intermediate values.  Therefore, with the exception of the case when
$\star = /$ and $0 \in Y$, the range $X \star Y$ has an interval form
$[\min{(x \star y)}, \max{(x \star y)}]$, where the $\min$ and $\max$ are taken over all
pairs $(x,y) \in X \times Y$.  Fortunately, we do not have to evaluate $x \star y$ over 
every $(x, y) \in X \times Y$ to find the global $\min$ and global $\max$ of $\star(x,y)$ 
over $X \times Y$, because the monotonicity of the $\star(x,y^*)$ in terms of $x \in X$ for 
any fixed $y^* \in Y$ implies that the extremal values are attained on the boundary of 
$X \times Y$, i.e., the set $\{\underline{x}, \underline{y}, \overline{x}$, and $\overline{y} \}$.  
Thus the theorem can be verified by examining the finitely many boundary cases. $\square$

An extremely useful property of interval arithmetic that is a direct consequence of 
Definition \ref{Df:intarith} is summarized by the following theorem.  
\begin{theorem}[Fundamental property of interval arithmetic]\label{Th:FundPropIntArith}
If $X \subseteq X^{\prime}$ and $Y \subseteq Y^{\prime}$ and $\star \in \{+,-,\cdot, / \}$, then
\[
X \star Y \subseteq X^{\prime} \star Y^{\prime}, 
\]
where we require that $0 \notin Y^{\prime}$ when $\star = /$.
\end{theorem}
{\bf Proof}:
\[
X \star Y = \{ x \star y: x \in X, y \in Y \} \subseteq \{ x \star y: x \in X^{\prime}, y \in Y^{\prime} \} = 
X^{\prime} \star Y^{\prime}. \square
\]
Note that an immediate implication of Theorem \ref{Th:FundPropIntArith} is that when $X = x$ and $Y = y$ are thin
intervals (real numbers $x$ and $y$), then $X^{\prime} \star Y^{\prime}$ will contain the result of the real arithmetic 
operation $x \star y$.
\begin{definition}[Range]
Consider a real-valued function $f: D \rightarrow \mathbb{R}$ where the domain $D \subseteq \mathbb{R}^n$.  
The range of $f$ over any $E \subseteq D$ is represented by $Rng(f;E)$ and defined to be the set
\[
Rng(f;E) \triangleq \{ f(x) : x \in E \}
\]
However, when the range of $f$ over any $X \in \mathbb{IR}^n$ such that $X \subseteq D$ is of interest, 
we will use the short-hand $f(X)$ for $Rng(f;X)$.
\end{definition}

\begin{definition}[Interval extension of subsets of $\mathbb{R}^n$] \label{D:Boxes}
For any Euclidean subset $\mathbf{\Theta} \subseteq \mathbb{R}^n$ let us denote its interval extension 
by $\mathbb{I}\mathbf{\Theta}$ and define it to be the set
\[
\mathbb{I}\mathbf{\Theta} \triangleq \{ X \in \mathbb{IR}^n : \underline{x}, \overline{x} \in \mathbf{\Theta} \}
\]
We refer the the $k$th interval of interval vector or box $X \in \mathbb{IR}^n$ by $X_k$.
\end{definition}

\begin{definition}[Inclusion isotony]
An box-valued map $F: D \rightarrow \mathbb{IR}^m$, where $D \in \mathbb{IR}^n$, is inclusion isotonic if it 
satisfies the property 
\[
\forall \, X \subseteq Y \subseteq D \implies F (X) \subseteq F(Y).
\]
\end{definition}

\begin{definition}[The natural interval extension]\label{D:NIE}
Consider a real-valued function $f: D \rightarrow \mathbb{R}$ given by a formula, where the 
domain $D \in \mathbb{IR}^n$.  If real constants, variables, and operations in $f$ are replaced 
by their interval counterparts, then one obtains
\[
F(X): \mathbb{I}D \rightarrow \mathbb{IR}.
\]
$F$ is known as the natural interval extension of $f$.  This extension is well-defined if we do not run into
division by zero.
\end{definition}

\begin{theorem}[Inclusion isotony of rational functions]\label{Thm:Rational}
Consider the rational function $f(x) = p(x)/q(x)$, where $p$ and $q$ are polynomials.  Let $F$ be its 
natural interval extension such that $F(Y)$ is well-defined for some $Y \in \mathbb{IR}$ and let
$X, X^{\prime} \in \mathbb{IR}$.  Then we have
\[
\begin{array}{ccl}
(i) & \text{\em Inclusion isotony:} & \forall \, X \subseteq X^{\prime} \subseteq Y \implies F(X) \subseteq F(X^{\prime}) 
\, \text{\em , and } \\ 
(ii) & \text{\em Range enclosure:}   & \forall \, X \subseteq Y \implies Rng(f; X) = f(X) \subseteq F(X).
\end{array}
\]
\end{theorem}
{\bf Proof} (cf.~\cite{Tucker2004}):
Since $F(Y)$ is well-defined, we will not run into division by zero, and therefore (i) follows from the repeated 
invocation of Theorem \ref{Th:FundPropIntArith}.  We can prove (ii) by contradiction.  Suppose $Rng(f; X) \nsubseteq F(X)$.  
Then there exists $x \in X$, such that $f(x) \in Rng(f;X)$ but $f(x) \notin F(X)$.  This in turn implies that 
$f(x)=F([x,x]) \notin F(X)$, which contradicts (i).  Therefore, our supposition cannot be true and we have 
proved (ii) $Rng(f; X) \subseteq F(X)$. $\square$  

\begin{definition}[Standard functions]
Piece-wise monotone functions, including exponential, logarithm, rational power, absolute value, and trigonometric functions,
constitute the set of standard functions
\[
\mathfrak{S} = \{ \, a^x, {log}_b(x), x^{p/q}, |x|, \sin(x), \cos(x), \tan(x), \sinh(x), \ldots, \arcsin(x), \ldots \, \}.
\]
\end{definition}
Such functions have well-defined interval extensions that satisfy inclusion isotony and {\em exact range enclosure}, i.e., 
$Rng(f;X) = f(X) = F(X)$.  Consider the following definitions for the interval extensions for some monotone functions 
in $\mathfrak{S}$ with $X \in \mathbb{IR}$,
\[
\begin{array}{lclr}
\exp(X)		&=& [\exp(\underline{x}), \exp(\overline{x})] 	& \\
\arctan(X)	&=& [\arctan(\underline{x}), \arctan(\overline{x})] 	& \\
\sqrt{(X)}	&=& [\sqrt{(\underline{x})}, \sqrt{(\overline{x})}] & \text{ if } 0 \leq \underline{x} \\
\log(X)		&=& [\log(\underline{x}), \log(\overline{x})]   & \text{ if } 0 < \underline{x} \\
\end{array}
\] 
and a piece-wise monotone function in $\mathfrak{S}$ with $\mathbb{Z}^+$ and $\mathbb{Z}^-$ representing the 
set of positive and negative integers, respectively.
\begin{equation}
X^n = 
\begin{cases}
[\underline{x}^n,\overline{x}^n] & \text{:   if } n \in \mathbb{Z}^+ \text{ is odd}, \notag \\
[{\langle X \rangle}^n,{|X|}^n] & \text{:   if } n \in \mathbb{Z}^+ \text{ is even}, \notag \\
[1,1] & \text{:   if } n = 0, \notag \\
[1/ \overline{x},1/ \underline{x}]^{-n} & \text{:   if } n \in \mathbb{Z}^- ; 0 \notin X \notag 
\end{cases}
\end{equation}

\begin{definition}[Elementary functions]\label{D:ElemFunc}
A real-valued function that can be expressed as a finite combination of constants, variables, arithmetic operations, 
standard functions and compositions is called an elementary function.  The set of all such elementary functions is
referred to as $\mathfrak{E}$.
\end{definition}

\begin{definition}[Directed acyclic graph (DAG) of a function]\label{D:DAG}
One can think of the process by which an elementary function $f$ is computed as the result of
a sequence of recursive operations with the subexpressions $f_i$ of $f$ where, $i=1,\ldots,n < \infty$.  This 
involves the evaluation of the subexpression $f_i$ at node $i$ with operands $s_{i_i},s_{i_2}$ from 
the sub-terminal nodes of $i$ given by the directed acyclic graph (DAG) for $f$
\begin{equation}\label{E:odotf}
s_i = \odot f_i \triangleq 
\begin{cases}
f_i(s_{i_1},s_{i_2}) & \text{:  if node $i$ has 2 sub-terminal nodes $s_{i_1},s_{i_2}$}  \\ 
f_i(s_{i_1}) & \text{:  if node $i$ has 1 sub-terminal node $s_{i_1}$}   \\ 
I(s_{i}) & \text{:  if node $i$ is a leaf or terminal node, } $I(x) = x$.   
\end{cases}
\end{equation}
The leaf or terminal node of the DAG is a constant or 
a variable and thus the $f_i$ for a leaf $i$ is set equal to the respective constant or variable.  The 
recursion starts at the leaves and terminates at the root of the DAG.  The DAG for an 
elementary $f$ with $n$ sub-expressions $f_1,f_2,\ldots,f_n$ is :
{\large 
\begin{equation}\label{E:DAGRec}
\{ \odot f_i \}_{i=1}^n \quad \rightarrowtail \quad \odot f_n = f(x), 
\end{equation}
}
where each $\odot f_i$ is computed according to \eqref{E:odotf}.
\end{definition}

\begin{figure}
 \makebox{\centerline{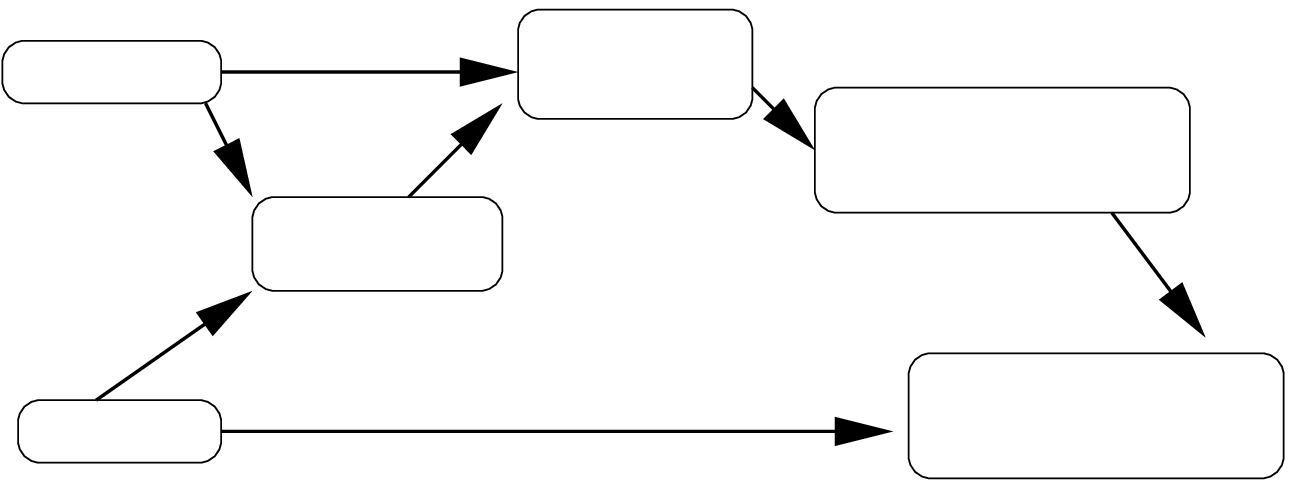}}
   \caption[]{Recursive evaluation of the sub-expressions $f_1,\ldots,f_6$ on the DAG
of the elementary function $f(x) = \odot f_6 =  x \cdot \sin ((x-3)/3)$}\label{Fi:DAG}
\end{figure}

For example the elementary function $x \cdot \sin ((x-3)/3)$
can be obtained from the terminus $\odot f_6$ of the recursion $\{ \odot f_i \}_{i=1}^6$ on the DAG for $f$ as shown 
in Fig.~\ref{Fi:DAG}.  
It would be convenient if guaranteed enclosures of the range $f(X)$ of an elementary $f$ can be 
obtained by its natural interval extension $F(X)$.  We show that inclusion
isotony does indeed hold for $F$, i.e.~if $X \subseteq Y$, then $F(X) \subseteq F(Y)$, and in
particular, the {\em inclusion property} that $x \in X \implies f(x) \in F(X)$ does hold.  

\begin{theorem}[The fundamental theorem of interval analysis]\label{3.1.11}
Consider any elementary function $f \in \mathfrak{E}$.  Let $F: Y \rightarrow \mathbb{IR}$ be its 
natural interval extension such that $F(Y)$ is well-defined for some $Y \in \mathbb{IR}$ and let 
$X, X^{\prime} \in \mathbb{IR}$.  Then we have
\[
\begin{array}{ccl}
(i) & \text{\em Inclusion isotony:} & \forall \, X \subseteq X^{\prime} \subseteq Y \implies F(X) \subseteq F(X^{\prime}) 
\, \text{\em , and } \\ 
(ii) & \text{\em Range enclosure:}   & \forall \, X \subseteq Y \implies Rng(f; X) = f(X) \subseteq F(X).
\end{array}
\]
\end{theorem}
{\bf Proof} (cf.~\cite{Tucker2004}): 
Any elementary function $f \in \mathfrak{E}$ is defined by the recursion \ref{E:DAGRec} on its 
sub-expressions $f_i$ where $i \in \{ 1,\ldots,n \}$ according to its DAG.
If $f(x) = p(x)/q(x)$ is a rational function, then the theorem already holds by Theorem \ref{Thm:Rational}, 
and if $f \in \mathfrak{S}$ then the theorem holds because the range 
enclosure is exact for standard functions.  Thus it suffices to show that if the theorem 
holds for $f_1, f_2 \in \mathfrak{E}$, then the theorem also holds for $f_1 \star f_2$, where 
$\star \in  \{ +,-,/,\cdot,\circ \}$.  By $\circ$ we mean the composition operator.  Since 
the proof is analogous for all five operators, we only focus on the $\circ$ operator.
Since $F$ is well-defined on its domain $Y$, neither the real-valued $f$ nor any of its 
sub-expressions $f_i$ have singularities in its respective domain $Y_i$ induced by $Y$.  
In particular $f_2$ is continuous on any $X_2$ and $X^{\prime}_2$ such that 
$X_2 \subseteq X^{\prime}_2 \subseteq Y_2$ implying the compactness of $F_2(X_2) \triangleq W_2$ and 
$F_2(X_2^{\prime}) \triangleq W_2^{\prime}$, respectively.  By our assumption that $F_1$ and $F_2$ 
are inclusion isotonic we have that $W_2 \subseteq W_2^{\prime}$ and also that 
\[
F_1 \circ F_2 (X_2) = F_1(F_2(X_2)) = F_1(W_2) \subseteq F_1(W_2^{\prime}) = F_1(F_2(X_2^{\prime})) = F_1 \circ F_2 (X_2)
\] 
The range enclosure is a consequence of inclusion isotony by an argument identical to that given in
the proof for Theorem \ref{Thm:Rational}. $\square$ 

The fundamental implication of the above theorem is that it allows us to enclose the range of
any elementary function and thereby produces an upper bound for the global maximum and a
lower bound for the global minimum over any compact subset of the domain upon which the function 
is well-defined.  We will see in the sequel that this is the work-horse of randomized 
enclosure algorithms that efficiently produce samples even from highly multi-modal target distributions.

\begin{figure}
\centering   \makebox{\includegraphics[width=5.0in]{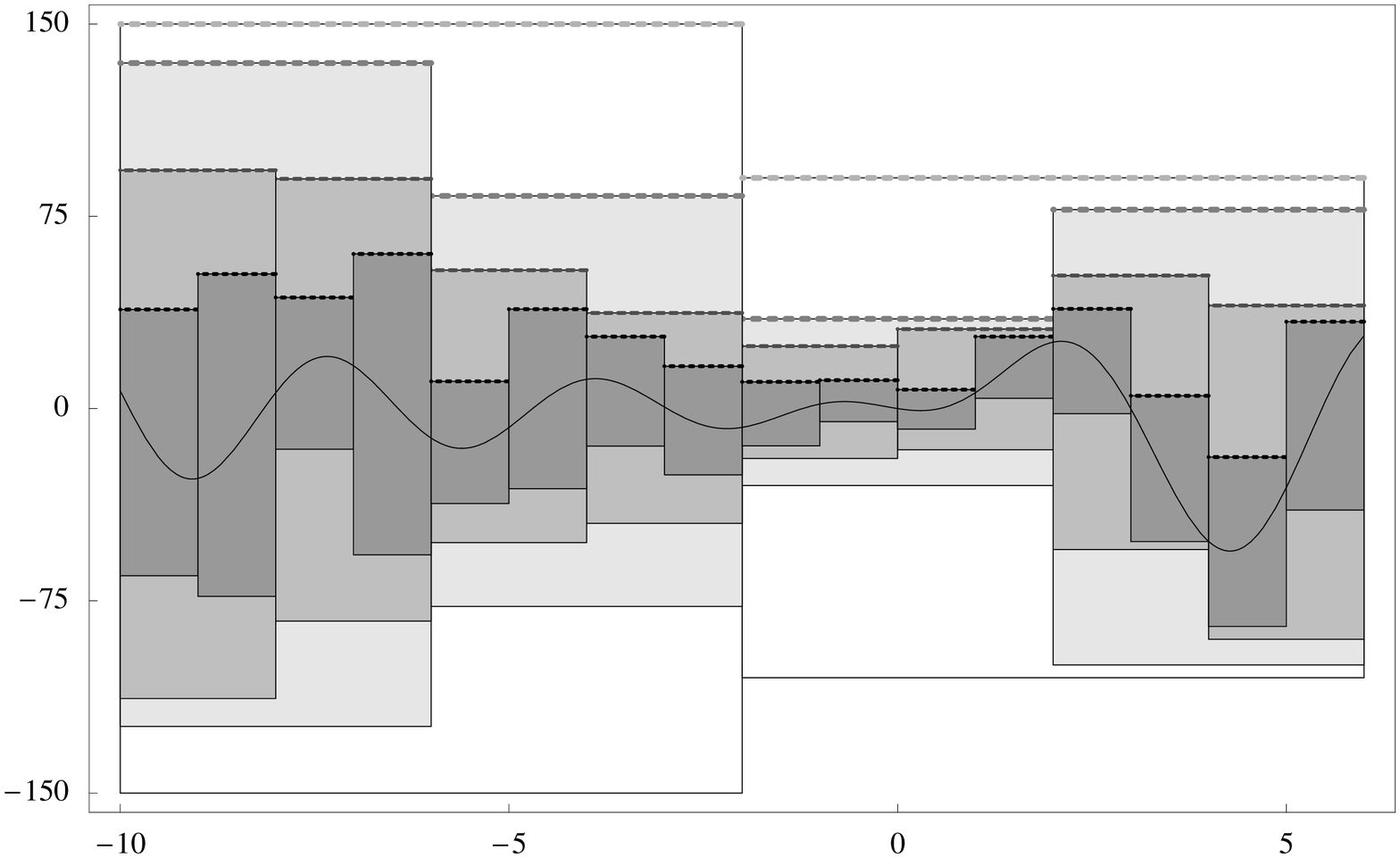}}  
   \caption[]{Range enclosure of the interval extension of 
   $-\sum_{k=1}^5{ k \, x \, \sin{(\frac{k(x-3)}{3})}}$ linearly tightens with the mesh.}\label{Fi:refine}
\end{figure}

Unlike the natural interval extension of an $f \in \mathfrak{S}$ that produces exact 
range enclosures, the natural interval extension $F(X)$ of an $f \in \mathfrak{E}$ often 
overestimates the range $f(X)$, but can be shown under mild conditions
to linearly approach the range as the maximal diameter of the box $X$ goes to zero, i.e.,
$\mathfrak{h}(F(X),f(X)) \leq \alpha \cdot d_{\infty}(X) \triangleq \max_{i}{d(X_i)}$ for some $\alpha \geq 0$.  This
implies that a partition of $X$ into smaller boxes $\{X^{(1)}, \cdots , X^{(m)}\}$ gives better 
enclosures of $f(X)$ through the union
$\bigcup_{i=1}^m F(X^{(i)})$ as illustrated in Fig.~\ref{Fi:refine}.  
Next we make the above statements precise.

\begin{definition}\label{D:LipElemFunc}
A function $f:D \rightarrow \mathbb{R}$ is Lipschitz if there exists
a Lipschitz constant $K$ such that, for all $x,y \in D$, we have $|f(x)-f(y)| \leq K |x-y|$.
We define $\mathfrak{E_L}$ to be the set of elementary functions whose sub-expressions $f_i$, $i=1,\ldots,n$ 
at the nodes of the corresponding DAGs are all Lipschitz.  
\end{definition}

\begin{theorem}[Range enclosure tightens linearly with mesh]\label{3.1.15}
Consider a function $f:D \rightarrow \mathbb{R}$ with $f \in \mathfrak{E_L}$.  Let $F$ be an inclusion isotonic interval extension
of $f$ such that $F(X)$ is well-defined for some $X \in \mathbb{IR}, X \subseteq I$.  Then there exists a 
positive real number $K$, depending on $F$ and $X$, such that if $X = \cup_{i=1}^k X^{(i)}$, then
\[
Rng(f;X) \subseteq \bigcup_{i=1}^k F(X^{(i)}) \subseteq F(X)
\]
and
\[
r \left( \bigcup_{i=1}^k F(X^{(i)}) \right) \leq r(Rng(f;X)) + K \max_{i=1,\ldots,k} {r(X^{(i)})}
\]
\end{theorem}
{\bf Proof }:  The proof is given by an induction on the DAG for $f$ similar to the proof of 
Theorem \ref{3.1.11} (See \cite{Tucker2004}).

\section{Appendix B}\label{S:AppB}
Here we will study the Moore rejection sampler (MRS) carefully.  Lemma \ref{lemma1} 
shows that MRS indeed produces independent samples from the desired target and 
Lemma \ref{lemma2} describes the asymptotics of the acceptance probability as the partition
of the domain is refined.

\begin{lemma}\label{lemma1}
Suppose that the target shape $p^*$ has a well-defined natural interval extension $P^*$.
If $U$ is generated according to Algorithm \ref{A:RS}, and if the proposal density $q^{\mathfrak{T}}(\theta)$  and the envelope function $f_{q^{\mathfrak{T}}}(\theta)$ 
are given by \eqref{E:qRS} and \eqref{E:fqRS}, respectively, then $U$ is distributed according to the target $p$.
\end{lemma}
{\bf Proof}:
From  \eqref{E:qRS} and \eqref{E:fqRS} observe that $f_{q^{\mathfrak{T}}}(t) = q^{\mathfrak{T}}(t) N_{q^{\mathfrak{T}}}$.
Let us define the following two subsets of $\mathbb{R}^2$,
\[
\mathcal{B}_q = \{(t,h) : 0 \leq h \leq f_{q^{\mathfrak{T}}}(t) \}, \text{   and   } 
\mathcal{B}_p = \{(t,h) : 0 \leq h \leq p^*(t) \}.
\]
First let us agree that Algorithm \ref{A:RS} produces a pair $(T,H)$ that is uniformly distributed on $\mathcal{B}_q$.  We can see this by
letting $k(t,h)$ denote the joint density of $(T,H)$ and $k(h | t)$ denote the conditional
density of $H$ given $T=t$.  Then,
\[
k(t,h) = 
\begin{cases}
q^{\mathfrak{T}}(t) \, k(h | t) & \text{ if } (t,h) \in \mathcal{B}_q \\
0   & \text{ otherwise }.
\end{cases}
\]
Since we sample a uniform height $h$ for a given $t$,
\[
k(h | t) = 
\begin{cases}
(f_{q^{\mathfrak{T}}}(t))^{-1} = (q^{\mathfrak{T}}(t)  N_{q^{\mathfrak{T}}})^{-1} 
& \text{ if } h \in [0, f_{q^{\mathfrak{T}}}(t)] \\
0  & \text{ otherwise}.
\end{cases}
\] 
Therefore,
\[
k(t,h) = 
\begin{cases}
q^{\mathfrak{T}}(t) \, k(h | t) = q^{\mathfrak{T}}(t) / (q^{\mathfrak{T}}(t) \, N_{q^{\mathfrak{T}}} )
= (N_{q^{\mathfrak{T}}})^{-1} & \text{ if } (t,h) \in \mathcal{B}_q \\
0   & \text{ otherwise }.
\end{cases}
\]
Thus we have shown that the joint density of $(T,H)$ is a uniformly distribution on $\mathcal{B}_q$.  The 
above relationship also makes geometric sense since the volume of $\mathcal{B}_q$ is exactly $N_{q^{\mathfrak{T}}}$.
Now, let $(T^*,H^*)$ be an accepted point, i.e., $(T^*,H^*) \in \mathcal{B}_p \subseteq \mathcal{B}_q$.  Then, the
uniform distribution of $(T,H)$ on $\mathcal{B}_q$ implies the uniform distribution of $(T^*,H^*)$ on $\mathcal{B}_p$.  
Since the volume of $\mathcal{B}_p$ is $N_p$, the p.d.f.~of $(T^*,H^*)$ is identically $1/N_p$ on $\mathcal{B}_p$ and 
$0$ elsewhere.  Hence, the marginal p.d.f. of $U=T^*$ is
\[
\begin{array}{lcl}
w(u) & = & \int_0^{p^*(u)} 1/N_p \, d h \\ 
     & = & 1/N_p \int_0^{p^*(u)} 1 \, d h \\
     & = & 1/N_p \int_0^{N_p p(u)} 1 \, d h, \quad \because \, p(u) = p^*(u)/N_p \\
     & = & p(u). \quad \square  
\end{array}
\]
  
\begin{lemma} \label{lemma2}
Let $\mathfrak{U}_W$ be the uniform partition of $\mathbf{\Theta} = [\underline{\theta},\overline{\theta}]$ 
into $W$ intervals each of diameter $w$
\[
\begin{array}{lcl}
w & = &  \frac{(\overline{\theta} - \underline{\theta})}{W} \\ 
\Theta_W^{(i)} & = & [ \ \underline{\theta} + (i-1) w, \ \underline{\theta} + i w \ ] \, , i = 1, \dots, W \\ 
\mathfrak{U}_W & = & \{ \Theta_W^{(i)} \, , i = 1, \dots, W \}. 
\end{array}
\]
and let $p^* \in \mathfrak{E_L}$, then
\[
\mathbf{A}^p_{\mathfrak{U}_W} = 1 - \mathcal{O} (1/W)
\]
\end{lemma}
{\bf Proof} \\
Then by means of Theorem \ref{3.1.15} 
\[
\begin{array}{lcl}
d(\Theta_W^{(i)}) = \mathcal{O} (1/W) & \implies & \mathfrak{h}( \ p^*(\Theta_W^{(i)}), P^*(\Theta_W^{(i)}) \ ) 
= \mathcal{O} (1/W) \\
& \implies & d(P^*(\Theta_W^{(i)}) )=  \mathcal{O} (1/W), \qquad \because p^* \in {\mathfrak{E_L}}
\end{array}
\]
Therefore
\[
  \sum_{i=1}^{|\mathfrak{U}_W|} \left( d(\Theta_W^{(i)}) \cdot P^*(\Theta_W^{(i)}) \right) 
= w \sum_{i=1}^{W} P^* \left( [ \ \underline{\theta} + (i-1) w, \ \underline{\theta} + i w \ ] \right),
\]
and we have 
\[
\begin{array}{lcl}
d(w \sum_{i=1}^{W} P^*(\Theta_W^{i})) = \mathcal{O} (1/W)
& \implies & \mathbf{A}^p_{\mathfrak{U}_W} = 1 - \mathcal{O} (1/W)
\end{array}
\]
Therefore the lower bound for the acceptance probability $\mathbf{A}^p_{\mathfrak{U}_W}$ 
of MRS approaches $1$ no slower than linearly with the refinement of 
$\mathbf{\Theta}$ by $\mathfrak{U}_W$.  
Note that this should hold for a general nonuniform partition with $w$ 
replaced by the mesh.$\square$
 
\bibliographystyle{splncs}
\bibliography{references} 
\end{document}